\documentclass [12pt]{article}
\usepackage {graphicx}
\usepackage{booktabs}
\usepackage{multirow}

\begin{document}

\title{A Short Tale of Long Tail Integration}

\author{Xiaolin Luo \\
\footnotesize{CSIRO Mathematical and
Information Sciences, Sydney,} \\
\footnotesize{Locked Bag 17, North Ryde, NSW, 1670, Australia.}\\
\footnotesize{e-mail: Xiaolin.Luo@csiro.au} \\
{}\\
Pavel V.~Shevchenko \\
\footnotesize{CSIRO Mathematical and
Information Sciences, Sydney,} \\ \footnotesize{Locked Bag 17, North
Ryde, NSW, 1670, Australia.}\\ \footnotesize{e-mail:
Pavel.Shevchenko@csiro.au}}

\date{First version 2 July 2009\\
Revised 5 March 2010}

\maketitle

\begin{abstract}
\noindent Integration of the form $\int_a^\infty {f(x)w(x)dx} $,
where $w(x)$ is either $\sin (\omega {\kern 1pt} x)$ or $\cos
(\omega {\kern 1pt} x)$, is widely encountered in many engineering
and scientific applications, such as those involving Fourier or
Laplace transforms. Often such integrals are approximated by a
numerical integration over a finite domain $(a,\,b)$, leaving a
truncation error equal to the tail integration $\int_b^\infty
{f(x)w(x)dx} $ in addition to the discretization error. This paper
describes a very simple, perhaps the simplest, end-point correction
to approximate the tail integration, which significantly reduces the
truncation error and thus increases the overall accuracy of the
numerical integration, with virtually no extra computational effort.
Higher order correction terms and error estimates for the end-point
correction formula are also derived. The effectiveness of this
one-point correction formula is demonstrated through several
examples.

\vspace{1cm} \noindent \textbf{Keywords:} numerical integration,
Fourier transform, Laplace transform, truncation error.
\end{abstract}

\pagebreak

\section{Introduction}
\label{sec:introductionords} Integration of the form $\int_a^\infty
{f(x)w(x)dx} $, where $w(x)$ is either $\sin (\omega {\kern 1pt} x)$
or $\cos (\omega {\kern 1pt} x)$, is widely encountered in many
engineering and scientific applications, such as those involving
Fourier or Laplace transforms. Often such integrals are approximated
by  numerical integrations over a finite domain $(a,\,b)$, resulting
in a truncation error  $\int_b^\infty {f(x)w(x)dx} $, in addition to
the discretization error. One example is a discrete  Fourier
transform (DFT), where there is a truncation error due to cut-off in
the tail, in addition to the discretization error.

In theory the cut-off error can always be reduced by extending the
finite domain at the expense of computing time. However, in many
cases a sufficiently long integration domain covering a very long
tail can be computationally expensive, such as when the integrand
$f(x)$ itself is a semi-infinite integration (e.g.  forward Fourier
or Laplace transform), or when the integrand decays to zero very
slowly (e.g.  a heavy tailed density or its characteristic
function). Much work has been done to directly compute the tail
integration in order to reduce the truncation error. Examples
include nonlinear transformation and extrapolation (Wynn
1956\nocite{Wynn56}, Alaylioglu et al
 1973\nocite{Alay73},   Sidi 1980\nocite{Sidi80}, 1982\nocite{Sidi82}, 1988\nocite{Sidi88}, Levin and Sidi 1981\nocite{LeviSi81}) and
application of special or generalized quadratures (Longman
1956\nocite{Long56},
 Hurwitz and Zweifel 1956\nocite{HurwZw56},
Bakhvalov and Vasileva 1968\nocite{BakhVas68}, Piessens
1970\nocite{Pies70},
 Piessens and Haegemans 1973\nocite{PiesHa73},
 Patterson 1976\nocite{Patt76},  Evans and Webster
1997\nocite{EvanWeb97},  Evans and Chung 2007\nocite{EvanChu07}),
among many others. This paper describes a very simple, perhaps the
simplest, end-point correction to account for the tail integration
over the entire range $(b,\infty )$. The treatment of the tail
reduces the usual truncation error significantly to a much smaller
discrete error, thus increasing overall accuracy of the integration,
while requiring virtually no extra computing effort. For the same
accuracy, this simple tail correction allows a much shorter finite
integration domain than would be required otherwise, thus saving
computer time while avoiding extra programming effort. To our
knowledge this result is not known in the literature and we believe
it deserves to be published for its elegant simplicity and broad
applicability. Though it is possible that our formula is a
rediscovery of a very old result hidden in the vast literature
related to numerical integration.

The paper is organized as follows. In Section 2, we derive the tail
integration approximation and its analytical error. A few examples
are shown to demonstrate the effectiveness of the tail integration
approximation in Section 3. Concluding remarks are given in Section
4.

\section{Tail integration}
Consider integration $\int_a^\infty {f(x)\sin (\omega {\kern 1pt}
x)dx} $. Without loss of generality, we assume $\omega = 1$ (a
change of variable $y = \omega {\kern 1pt} x$ results in the desired
form). For $\int_a^\infty {f(x)\cos (\omega {\kern 1pt} x)dx} $ the
derivation procedure and the resulting formula are very similar. In
the following, we assume that

\begin{itemize}
                         \item The integral $\int_a^\infty {f(x)\sin (\omega {\kern 1pt}
x)dx} $ exists;
                         \item All derivatives $f^{(k)}(x)$ exist
                         and $\to 0$ as $k \to \infty$.
                       \end{itemize}

\subsection{Piecewise linear approximation}
The truncation error of replacing $\int_a^\infty {f(x)\sin (x)dx} $
by $\int_a^b {f(x)\sin (x)dx} $ is simply the tail integration

\begin{equation}
I_T = \int\limits_b^\infty {f(x)\sin (x)dx}.
\end{equation}

\noindent For higher accuracy, instead of increasing truncation
length at the cost of computing time, we propose to compute the tail
integration $I_T $ explicitly by a very economical but effective
simplification. Assume $f(x)$ approaches zero as $x \to \infty $ and
the truncation point $b$ can be arbitrarily chosen in a numerical
integration. Let $b = N\pi $, where $N$ is some large integer.
Dividing integration from $N\pi $ to $\infty $ into cycles with an
equal length of $\pi $ yields

\begin{equation}
\int\limits_{N\pi }^\infty {f(x)\sin (x)dx} = \sum\limits_{k =
0}^\infty {I_k {\kern 1pt} ,\quad I_k = \int\limits_{(N + k)\pi
}^{(N + k + 1)\pi } {f(x)\sin (x)dx} }.
\end{equation}

\noindent Now assume that $f(x)$ is piecewise linear within each
$\pi $-cycle, so that each of the integrals $I_k $ in (2) can be
computed exactly. That is, in the range $\left[ {(N + k)\pi , (N + k
+ 1)\pi } \right]$, we assume that $f(x)$ is approximated by

\begin{equation}
f(x) \approx f_k + \frac{x - (N + k)\pi }{\pi }(f_{k + 1} - f_k ),
\end{equation}

\noindent where   $f_k=f((N+k)\pi)$. Substitute (3) into (2), then
analytical integration by parts of each $I_k$ in (2) gives

\begin{eqnarray}
\int\limits_{N\pi }^\infty f(x)\sin (x)dx =\sum\limits_{k =0}^\infty
{I_k}  \approx  \sum\limits_{k = 0}^\infty {( - 1)^{N + k}(f_k  +
f_{k + 1} )} = ( - 1)^Nf_0 = ( - 1)^Nf(N\pi ).
\end{eqnarray}

\noindent This elegant result given by (4) means that we only need
to evaluate the integrand $f(x)$ at one single point $x = N\pi $
(the truncation point) for the entire tail integration, replacing
the truncation error with a much smaller round-off error. As will be
demonstrated later, this one-point formula for the potentially
demanding tail integration is remarkably effective in reducing the
truncation error caused by ignoring $I_T $.

\subsection{Higher order correction terms and error estimates}
Formula (4) can be derived more generally through integration by
parts, and a recursive deduction gives us higher order correction
terms and thus error estimates. Integrating (1) by parts with $b =
N\pi $, we have

\begin{equation}\int\limits_{N\pi }^\infty {f(x)\sin (x)dx = \,} ( - 1)^Nf(N\pi
) + \int\limits_{N\pi }^\infty {f^{(1)}(x)\cos (x)dx},
\end{equation}

\noindent where $f^{(1)}(x) = df(x) / dx$. If we assume $f(x)$ is
linear within each $\pi $-cycle in the tail, then the integration
$\int_{N\pi }^\infty {f^{(1)}(x)\cos (x)dx} $ vanishes, because
within each $\pi $-cycle $f^{(1)}(x)$ is constant from the piecewise
linear assumption and $\int_{k\pi }^{(k + 1)\pi } {\cos (x)dx} = 0$
for any integer $k$, and $f^{(1)}(\infty ) \to 0$ as $f(\infty ) \to
0$. Thus, under the piecewise linear assumption, (5) and (4) are
identical. Continuing with integration by parts in (5) and noting
$f^{(1)}(x) \to 0$ at infinity, we further obtain

\begin{equation}
\int\limits_{N\pi }^\infty {f(x)\sin (x)dx = } \,( - 1)^Nf(N\pi ) -
\int\limits_{N\pi }^\infty {f^{(2)}(x)\sin (x)dx},
\end{equation}

\noindent where $f^{(2)}(x) = d^2f(x) / dx^2$. Equation (6), as well
as (5), is exact -- no approximation is involved. The recursive
pattern in (6) is evident. If we now assume that the second
derivative $f^{(2)}(x)$ is piecewise linear in each $\pi $-cycle in
the tail, then (6) becomes

\begin{equation}
\int\limits_{N\pi }^\infty {f(x)\sin (x)dx \approx \,} ( -
1)^N\left( {f(N\pi ) - f^{(2)}(N\pi )} \right).
\end{equation}

With the additional correction term, (7) is more accurate than (4).
In general, without making any approximation, from the recursive
pattern of (6) we arrive at the following expression for the tail
integral

\begin{equation}\int\limits_{N\pi }^\infty {f(x)\sin (x)dx} =
( - 1)^N f(N\pi ) + \sum\limits_{i = 1}^{k-1} {( - 1)^{N +
i}f^{(2i)}(N\pi )}+ (-1)^k\int\limits_{N\pi }^\infty
{f^{(2k)}(x)\sin (x)dx},
\end{equation}

\noindent where $k\geq 1$, $f^{(2k)}(N\pi)$ is the 2$k$-th order
derivative of $f(x)$ at the truncation point. As will be shown later
with examples, typically the first few terms from (8) are
sufficiently accurate. The error in using formula (4) is readily
obtained from (8)

\begin{equation}\varepsilon_T= \sum\limits_{i = 1}^{k-1} {( - 1)^{N +
i}f^{(2i)}(N\pi )}+ (-1)^k\int\limits_{N\pi }^\infty
{f^{(2k)}(x)\sin (x)dx},
\end{equation}

 In deriving (8), we have assumed all derivatives
exist and $f^{(k)}(\infty ) = 0$. Under certain conditions, the
infinite series in (8) and (9) represents the integral
asymptotically as $N \to \infty $, i.e. we have the asymptotic
expansion

$$\int\limits_{N\pi }^\infty {f(x)\sin (x)dx}\sim \sum\limits_{i = 0}^{\infty} {( - 1)^{N +
i}f^{(2i)}(N\pi )}\;\;\; \mbox{as} \;\;\; N \to \infty$$

For example, if we assume that, for some $\gamma<0$,

$$f^{(m)}(x)=O(x^{\gamma-m})\;\;\; \mbox{as} \;\;\;
N\to\infty,\;\;\;m=0,1,2,...,$$

\noindent then the integral term on the right-hand side of (8) can
be bounded by $C_kN^{\gamma-2k+1}$ as $N\to\infty$, for some
positive constant $C_k$, and the series converges to the integral.

The derivatives approaching zero as $x \to \infty $ is a consequence
of the existence of integral (1). Otherwise, if  $f^{(1)}(\infty )
\neq 0$,  integral (1) does not exist, which is evident form (5).
Applying this argument recursively, all derivatives $f^{(k)}(\infty
) = 0$, if they exist. Obviously if $f(x)$ is a power function (e.g.
$1/\sqrt{x}$), the ratio $f^{(1)}(x) / f(x)$ is of the order $1/x$
as $x \to \infty $, so is the ratio $f^{(k+1)}(x) / f^{(k)}(x)$.
This implies that, for a power-like function, each error term in (9)
decreases by two orders of magnitude from its preceding term as the
index number $k$ increases by one.

\vspace{0.4cm}
 \noindent\textbf{Remark}. Note that there is no
truncation error in (4) and the error is a discretization error in
nature. In theory, the tail integration error can be estimated by
(9). In practice, however, derivatives of integrand at the
truncation point may only be evaluated numerically. The assumption
of piecewise linearity, although reasonable for $f(x)$ at large $x$,
may appear to be rather crude for a high precision computation.
However, we recall that we are only trying to reduce the already
small truncation error $I_T $ and a reasonable approximation in $I_T
$ could lead to significant improvement in the overall accuracy of
integration. For example, suppose a relative error of 1{\%} due to
ignoring truncation and 10{\%} error in evaluating the tail
integration using the very simple formula (4). The overall accuracy
with this tail integration added is now improved from 1{\%} to
0.1{\%} (1{\%} times 10{\%}). This improvement by an order of
magnitude is achieved by simply evaluating the integrand at the
truncation point. The assumption of a piecewise linearity applies to
a broad range of functions, thus the special tail integration
approximation can have a wide application. Note, piecewise linear
assumption does not even require monotonicity - $f(x)$ can be
oscillating, as long as its frequency is relatively small compared
with the principal cycles in $\sin (x)$, as demonstrated in one of
the examples below.

If the oscillating factor is $\cos (x)$ instead of $\sin (x)$, we
can still derive a one-point formula similar to (4) by starting the
tail integration at $(N - 1 / 2)\pi $ instead of $N\pi $. In this
case, the tail integration is

\begin{equation}\int\limits_{(N - 1 / 2)\pi }^\infty {f(x)\cos (x)dx \approx }
( - 1)^Nf\left( {(N - 1 / 2)\pi } \right).
\end{equation}

\noindent Also, the tail integration approximation can be applied to
the left tail (integrating from $ - \infty $ to $ - N\pi )$ as well,
if such integration is required.

It is known from the literature that truncation is better at extrema
 of the oscillatory part than at the zeros (Lyness 1986,\nocite{Lyness86} Espelid and Overholt
1994 \nocite{Espelid94} and Sauter 2000\nocite{Sauter2000}).
Truncating at $x_T=(N-1/2)\pi$, the extrema for $sin(x)$, we obtain
an expression for the tail integration or the truncation error
similar to (8)

\begin{eqnarray}\int\limits_{x_T }^\infty {f(x)\sin (x)dx} &=&
( - 1)^N f^{(1)}(x_T ) + \sum\limits_{i = 1}^{k-1} {( - 1)^{N +
i}f^{(2i+1)}(x_T)} \nonumber \\
&& + (-1)^k\int\limits_{x_T}^\infty
{f^{(2k)}(x)\sin (x)dx}.
\end{eqnarray}

The leading term of the truncation error is now $f^{(1)}(x_T)$ in
(11), compared with $f(x_T+1/2\pi)$ in (8). Assuming
$|f^{(1)}(x_T)/f(x_T+1/2\pi)|<1$ for some large $x_T$, e.g. when
$f(x)$ is a power-like function, then it is obvious truncation at
extrema has a smaller truncation error. However, our formula is
about the reduction of the truncation error by including an
approximation of the tail integration. If truncation is done at
$x=(N-1/2)\pi$ instead of $x=N\pi$, then the first correction term
will be $(-1)^{N+1}f^{(1)}((N-1/2)\pi)$, involving the first
derivative of $f(x)$. In many important applications the first
derivative of $f(x)$ cannot be evaluated accurately. For example,
when inverting a characteristic function of a compound distribution,
$f(x)$ itself is a semi-infinite integration of an oscillatory
function, which could only be obtained numerically. Taking finite
difference of a numerically evaluated function will in general
reduce the accuracy by an order of magnitude. So for general
purposes the truncation is chosen at the zeros, i.e. at $x=N\pi$.

Of course, if derivative of $f(x)$  is in closed form and can be
accurately evaluated, truncation and correction at extrema will
indeed be more accurate, with a leading error term of
$f^{(3)}(x_T)$. But we could also include the second derivative term
for the truncation at zeros, with a leading error term of
$f^{(4)}(x_T)$, and so on. In general when higher order derivatives
can be computed precisely, then one can include some higher order
terms to reduce truncation error further and it does not matter much
whether the truncation is done at extrema or at zeros.

\section{Examples of tail integration}
The effectiveness of the above tail integration approximation is now
demonstrated in a few examples. Introduce the following notations

\begin{eqnarray}
I_E &=& \int_0^\infty {f(x)\sin (x)dx},\nonumber\\
\tilde {I}(N\pi ) &=& \int_0^{\pi {\kern 1pt} N} {f(x)\sin
(x)dx},\nonumber\\
I_T (N\pi ) &=& \int_{\pi {\kern 1pt} N}^\infty {f(x)\sin
(x)dx}.\nonumber
\end{eqnarray}

\noindent In all the following examples the exact semi-infinite
integration $I_E $ is known in closed form, and its truncated
counterpart $\tilde {I}(N\pi )$ is either known in closed form or
can be computed accurately. For simplicity in all the examples $N$
is taken to be an even number, i.e $(-1)^N=1$. The exact tail
integration $I_T (N\pi )$ can be computed from $I_T (N\pi ) = I_E -
\tilde {I}(N\pi )$. We compare $\tilde {I}(N\pi ) + f(N\pi )$ with
$\tilde {I}(N\pi )$ and compare both of them with the exact
semi-infinite integration $I_E $. The error reduction can be
quantified by comparing the ``magic'' point value given by formula
(4) with the exact tail integration $I_T (N\pi )$. Also note that
the analytic formula for the error of using (4), $\varepsilon _T =
I_E - [\tilde {I}(N\pi ) + f(N\pi )]$, is given by (9).

~

\noindent\textbf{\textit{Example 1: }}$f(x) = e^{ - \alpha
x},\,\;(\alpha
> 0).$

\noindent In this example, the closed form results are

$$
I_E = \int\limits_0^\infty {e^{ - \alpha x}\sin (x)dx = } \frac{1}{1
+ \alpha ^2},\quad \tilde {I}(N\pi ) = \int\limits_0^{N\pi } {e^{ -
\alpha x}\sin (x)dx = } \frac{1 - e^{ - \alpha \,N\pi }}{1 + \alpha
^2}.
$$

Figure 1 compares the ``magic'' point value $f(N\pi )$ representing
simplified tail integration with the exact tail integration $I_T
(N\pi ) = I_E - \tilde {I}(N\pi )$ as functions of parameter $\alpha
$ for $N = 10$, i.e. the truncated lengths $l_T = 10\pi $. The
figure shows that a simple formula (4) matches the exact
semi-infinite tail integration surprisingly well for the entire
range of parameter $\alpha $. Corresponding to Figure 1, the actual
errors of using formula (4) are shown in Table 1, in comparison with
the truncation errors without applying the correction term given by
(4). Figure 2 shows the same comparison at an even shorter truncated
length of $4\pi $. The error of using (4) is $\left|\varepsilon
_T\right| = {\alpha ^2\exp ( - \alpha \pi N) / (1 + \alpha ^2)} $.

If $\alpha $ is large, the function $f(x) = e^{ - \alpha x}$ is
``short tailed'' and it goes to zero very fast. The absolute error
$\left|\varepsilon_T\right| $ is very small even at $N = 4$. The
relative error (against the already very small tail integration),
given by $\left|\varepsilon_T\right| / \left| {I_E - \tilde
{I}(2N\pi )} \right| = \alpha ^2$, is actually large in this case.
But this large relative error in the tail approximation does not
affect the high accuracy of the approximation for the whole
integration. What is important is the error of the tail integration
relative to the whole integration value. Indeed, relative to the
exact integration, the error of using (4) is
$\left|\varepsilon_T\right| / \left| {I_E } \right| = \alpha ^2\exp
( - \alpha \pi N)$, which is about $2.7\times 10^{ - 53}$ at $N =
4$. The condition $f^{(k)}(x)\to 0$ as $k \to \infty$ is not
satisfied in this case if $\alpha > 1$. However, as discussed above,
the application of formula (4) does not cause any problem.

For a small value of parameter $\alpha $, the truncation error will
be large unless the truncated length is very long. For instance,
with $\alpha = 0.01$ the truncation error (if ignore the tail
integration) is more than 70{\%} at $l_T = 10\pi $ ($N = 10$, as the
case in Figure 1), and it is more than 88{\%} at $l_T = 4\pi $ ($N =
4$, as the case in Figure 2). On the other hand, if we add the
``magic'' value from formula (4) to approximate the tail
integration, the absolute error of the complete integration
$\left|\varepsilon_T\right|  $ due to this approximation is less
than 0.01{\%}, and the relative error is $\left|\varepsilon_T\right|
= \alpha ^2 = 0.01\% $ at both $l_T = 10\pi $ and $l_T = 4\pi $. In
other words, by including this one-point value, the accuracy of
integration has dramatically improved by several orders of magnitude
at virtually no extra cost, compared with the truncated integration.
For the truncated integration $\tilde {I}(N\pi )$ to have similar
accuracy as $\tilde {I}(4\pi ) + f(4\pi )$, we need to extend the
truncated length from $4\pi $ to $300\pi $ for this heavy tailed
integrand.

~

\noindent\textbf{\textit{Example 2: $f(x) = 1 / \sqrt x$ }}.

\noindent This example has a heavier tail than the previous one.
Here, we have closed form for $I_E $, but not for $\tilde {I}$ or
$I_T $,

$$
I_E = \int\limits_0^\infty {\frac{\sin (x)}{\sqrt x }dx = } \sqrt
{\frac{\pi }{2}} ,\quad \tilde {I}(N\pi ) = \int\limits_0^{N\pi }
{\frac{\sin (x)}{\sqrt x }dx}.
$$

\noindent $\tilde {I}$ or $I_T $ can be accurately computed by
adaptive integration functions available in many numerical packages.
Here we used \textit{IMSL} function based on the modified
Clenshaw-Curtis integration method (Clenshaw and Curtis
1960\nocite{ClCu60}; Piessens, Doncker-Kapenga, \"{U}berhuber and
Kahaner 1983\nocite{PiDoKaUbKa83}).

Figure 3 compares the ``exact'' tail integration $I_T (N\pi ) =
\int_{N\pi }^\infty {\sin (x) / \sqrt x dx} $ with the one-point
value $f(N\pi )$. Again the one-point approximation does an
extremely good job. Even at the shortest truncation length of just
$2\pi $ the one-point approximation is very close to the exact
semi-infinite tail integration. Applying the analytical error
formula (9) to $f(x) = 1 / \sqrt x $\textbf{, }we have

$$
\varepsilon _T = \sum\limits_{k = 1}^\infty {( - 1)^k\frac{1\times
3\times ...\times (4k - 1)}{2^{2k}x^{(4k + 1) / 2}}}.
$$

Taking the first three leading terms we get $\varepsilon _T \approx
- \mbox{0.00}6\mbox{95}$ at $N = 2$ and $\varepsilon _T \approx -
\mbox{2.392}\times \mbox{10}^{-5}$ at $N = 20$. The relative error
$\left|\varepsilon_T\right|  / I_E (2N\pi )$ is about 1{\%} at $N =
2$ and it is about 0.002{\%} at $N = 20$. Apparently, if the extra
correction term $f^{(2)}(N\pi )$ is included as in (7), the error
$\varepsilon_T$ reduces further by an order of magnitude at $N = 2$
and by several orders of magnitude at $N = 20$. Corresponding to
Figure 3, the actual errors of using formula (4) are shown in Table
2, in comparison with the truncation errors without applying the
correction term given by (4).

Figure 4 shows the truncated integration $\tilde {I}(N\pi )$ and the
truncated integration with the tail modification (4) added, i.e.
$\tilde {I} + f(N\pi )$, along with the correct value of the full
integration $I_E = \sqrt {\pi / 2} $. The contrast between results
with and without the one-point tail approximation is striking. At
the shortest truncation length of $2\pi $ ($N = 2)$, the relative
error due to truncation for the truncated integration $(I_E - \tilde
{I}(N\pi )) / I_E $ is more than 30{\%}, but with the tail
approximation added, the relative error $(I_E - \tilde {I}(N\pi ) -
f(N\pi )) / I_E $ reduces to 0.5{\%}. At $100\pi $, the largest
truncation length shown in Figure 4, the relative error due to
truncation is still more than 4{\%}, but after the ``magic'' point
value is added the relative error reduces to less than $0.5\times
10^{ - 6}$.

Another interesting way to look at these comparisons, which is relevant for
integrating heavy tailed functions, is to consider the required truncation
length for the truncated integration to achieve the same accuracy as the one
with the ``magic'' value added. For the truncated integration $\tilde
{I}(N\pi )$ to achieve the same accuracy of $\tilde {I}(2\pi ) + f(2\pi )$
(integration truncated at one-cycle plus the ``magic point value), we need
to extend the integration length to $7700\pi $. For $\tilde {I}(N\pi )$ to
achieve the same accuracy of $\tilde {I}(100\pi ) + f(100\pi )$, the
integration length has to be extended to more than $10^{12}\pi $! On the
other hand, if we add the tail approximation $f(7700\pi )$ to $\tilde
{I}(7700\pi )$, the relative error reduces from 0.5{\%} to less than $10^{ -
11}$! This error reduction requires no extra computing, since $f(7700\pi )$
is simply a number given by $1 / \sqrt {7700\pi } $.

~

\noindent\textbf{\textit{Example 3: }}$f(x) = \cos (\alpha x) /
x,\;\alpha < 1$\textbf{.}

\noindent We have remarked that the piecewise linear assumption does
not require monotonicity, i.e. $f(x)$ can be oscillating, as long as
its frequency is relatively small compared with the principal
cycles. For example, when the function $f(x)$ is the characteristic
function of a compound distribution, it oscillates with its
frequency approaching zero in the long tail. In the current example
with $f(x) = \cos (\alpha x) / x$, there is a closed form for $I_E
$, but not for $\tilde {I}$ or $I_T $,

$$
I_E = \int\limits_0^\infty {\frac{\cos (\alpha x)\sin (x)}{x}dx = }
\frac{\pi }{2},\quad \tilde {I}(2N\pi ) = \int\limits_0^{N\pi }
{\frac{\cos (\alpha x)\sin (x)}{x}dx} ,\quad \alpha < 1.
$$

Figure 5 compares the ``exact'' tail integration $I_T (N\pi )$ with
the one-point approximation $f(N\pi )$ for the case $\alpha = 0.2$.
Again the one-point approximation performs surprisingly well,
despite $f(x)$ itself is now an oscillating function, along with the
principal cycles in $\sin (x)$. The piecewise linearity assumption
is apparently still valid for relatively mild oscillating $f(x)$.
Corresponding to Figure 5, the actual errors of using formula (4)
are shown in Table 3, in comparison with the truncation errors
without applying the correction term given by (4). Not surprisingly,
the errors are larger in comparison with those in examples 1 and 2,
due to the fact that $f(x)$ now is itself an oscillating function.
Still, Table 3 shows the truncation error is reduced by an order of
magnitude after applying the simple formula (4).

 Figure 6 compares
the truncated integration $\tilde {I}(N\pi )$ against $\tilde {I} +
f(N\pi )$, along with the correct value of the full integration $I_E
= \pi / 2$. At truncation length $6\pi $, the shortest truncation
length shown in Figures 5 and 6, the relative error
$\left|\varepsilon_T\right| / I_E $ is less than 0.06{\%} and it is
less than 0.01{\%} at $N = 100$. In comparison, the truncated
integration without the end point correction has relative error of
2.7{\%} and 0.2{\%}, respectively for those two truncation lengths.
Applying the analytical error formula (9) to $f(x) = \cos (\alpha x)
/ x$ and noting $\sin (\alpha x) = 0$ and $\cos (\alpha x) = 1$ with
$\alpha = 0.2$ and $x = 100\pi $, we obtain

$$
\varepsilon_T \approx - \left( { - \frac{\alpha ^2}{x} +
\frac{2}{x^3}} \right) + \left( {\frac{\alpha ^4}{x} -
\frac{12\alpha ^2}{x^3} + \frac{24}{x^5}} \right),\quad x = l_T =
100\pi,
$$

\noindent where only the first two leading terms corresponding to
the 2$^{nd}$ and 4$^{th}$ derivatives are included, leading to
$\varepsilon_T \approx 0.0001273 + 5.07749\times 10^{- 6}$ at $N =
50$ that agrees with the actual error. Similar to the previous
example, if we include the extra correction term $f^{(2)}(N\pi )$,
the error reduces further by two orders of magnitude at $N = 100$.

The purpose of Example 3 is to show that the piecewise linear
approximation in the tail could still be valid even if there is a
secondary oscillation in $f(x)$, provided its frequency is not as
large as the principal oscillator. If the parameter $\alpha$ is
larger than one, then we can simply perform a change of variable
with $y=\alpha x$ and integrate $(1/y) \sin (y/\alpha) \cos (y)$ in
terms of $y$. Better still, for any value of $\alpha$, we can make
use of the equality $\cos (\alpha x)\sin (x)=(\sin(x-\alpha
x)+\sin(x+\alpha x))/2$ to get rid of the secondary oscillation
altogether before doing numerical integration. In practice, the
secondary oscillation often has a varying frequency with a slowly
decaying magnitude, such as in the case of the characteristic
function of a compound distribution with a heavy tail. In this case
it might be difficult to effectively apply regular numerical
quadratures in the tail integration, but the simple one-point
formula (4)  might be very effective.

\vspace{0.3cm}

 All these examples show dramatic
reduction in truncation errors if tail integration approximation (4)
is employed, with virtually no extra cost. If the extra correction
term $f^{(2)}(2N\pi )$ is included, i.e. using (7) instead of (4),
the error is reduced much further.

\section{Conclusions}
\label{sec:conclusionsxamples} We have derived perhaps the simplest
but efficient tail integration approximation, first intuitively by
piecewise linear approximation, then more generally through
integration by parts. Analytical higher-order correction terms and
thus error estimates are also derived. The usual truncation error
associated with a finite length of the truncated integration domain
can be reduced dramatically by employing the one-point tail
integration approximation, at virtually no extra computing cost, so
a higher accuracy is achieved with a shorter truncation length.

Under certain conditions outlined in the present study, the method
can be used in many practical applications. For example, the authors
have successfully applied the present method in computing heavy
tailed compound distributions through inverting their characteristic
functions, where the function $f(x)$ itself is a semi-infinite
numerical integration (Luo, Shevchenko and Donnelly
2007\nocite{LuShDo07}).

Of course there are more elaborate methods in the literature which
are superior to the present simple formula in terms of better
accuracy and broader applicability, such as some of the
extrapolation methods proposed by Wynn 1956 and by Sidi 1980, 1988.
The merit of the present proposal is its simplicity and
effectiveness - a single function evaluation for the integrand at
the truncation point  is all that is needed to reduce the truncation
error, often by orders of magnitude. It can not be simpler than
that. Also, in some applications the function $f(x)$ may not even
exist in closed form, for instance when $f(x)$ is the characteristic
function of some compound distributions as mentioned above, then
$f(x)$ itself is a semi-infinite integration of a highly oscillatory
function, which could only be obtained numerically. In such cases
some of the other more sophisticated methods relying on a closed
form of $f(x)$ may not be readily applicable.

\section{Acknowledgement}
We would like to thank David Gates, Mark Westcott and three
anonymous refrees for many constructive comments which have led to
significant improvements in the manuscript.


\begin{thebibliography}{10}

\bibitem{Alay73}
A.~Alaylioglu, G.~A. Evans, and J.~Hyslop, \emph{The evaluation of
oscillatory
  integrals with infinite limits}, J. Comp. Phys, \textbf{13} (1973), 433--438.

\bibitem{BakhVas68}
N.~S. Bakhvalov and L.~G. Vasil'eva, \emph{Evaluation of integrals
of
  oscillating functions by interpolation at nodes of gaussian quadratures},
  USSR Comp. Math. Math. Phys. \textbf{8} (1968), 241--249.

\bibitem{ClCu60}
C.~W. Clenshaw and A.~R. Curtis, \emph{A method for numerical
integration on an
  automatic computer}, Num. Math \textbf{2} (1960), 197--205.

\bibitem{Espelid94}
T.~O. Espelid and K.~J. Overholt, \emph{Dqainf: An algorithm for
automatic
  integration of infinite oscillating tails}, Numer. Algorithms \textbf{8}
  (1994), 83--101.

\bibitem{EvanChu07}
G.~A. Evans and K.~C. Chung, \emph{Evaluating infinite range
oscillating
  integrals using generalized quadrature methods}, Appl. Numer. Math.
  \textbf{57} (2007), 73--79.

\bibitem{EvanWeb97}
G.~A. Evans and J.~R. Webster, \emph{A high order progressive method
for the
  evaluation of irregular oscillating integrals}, Appl. Numer. Anal.
  \textbf{23} (1997), 205--218.

\bibitem{HurwZw56}
H.~Jr. Hurwitz and P.~F. Zweifel, \emph{Numerical quadrature of
fourier
  transform integrals}, MTAC \textbf{10} (1956), 140--149.

\bibitem{LeviSi81}
D.~Levin and A.~Sidi, \emph{Two new classes of nonlinear
transformations for
  accelerating the convergence of infinite integrals and series}, Appl. Math.
  Comput. \textbf{9} (1981), 175--215.

\bibitem{Long56}
I.~M. Longman, \emph{Note on a method for computing infinite
integrals of
  oscillatory functions}, Camb. Phil. Soc. Proc. \textbf{52} (1956), 764.

\bibitem{LuShDo07}
X.~Luo, P.~V. Shevchenko, and J.~Donnelly, \emph{Addressing impact
of
  truncation and parameter uncertainty on operational risk estimates}, The
  Journal of Operational Risk \textbf{2} (2007), 3--26.

\bibitem{Lyness86}
J.~Lyness and G.~Hines, \emph{To integrate some infinite oscillating
tails},
  ACM Trans. Math. software \textbf{12} (1986), 24--25.

\bibitem{Patt76}
T.~N.~L. Patterson, \emph{On high precision methods for the
evaluation of
  fourier integrals with finite and infinite limits}, Numer. Math. \textbf{27}
  (1976), 41--52.

\bibitem{Pies70}
R.~Piessens, \emph{Gaussian quadrature fomulas for the integration
of
  oscillating functions}, Math. Comp. \textbf{24} (1970), microfiche.

\bibitem{PiDoKaUbKa83}
R.~Piessens, E.~De. Doncker-Kapenga, C.~W. \"{U}berhuber, and D.~K.
Kahaner,
  \emph{Quadpack -- a subroutine package for automatic integration},
  Springer-Verlag, 1983.

\bibitem{PiesHa73}
R.~Piessens and A.~Haegemans, \emph{Numerical calculation of fourier
transform
  integrals}, Electron. Lett. \textbf{9} (1973), 108--109.

\bibitem{Sauter2000}
T.~Sauter, \emph{Computation of irregularly oscillating integrals},
Appl.
  Numer. Math. \textbf{35} (2000), 245--264.

\bibitem{Sidi80}
A.~Sidi, \emph{Extrapolation methods for oscillatory infinite
integrals}, J.
  Inst. Maths. Appl. \textbf{26} (1980), 1--20.

\bibitem{Sidi82}
\bysame, \emph{The numerical evaluation of very oscillatory
integrals by
  extrapolation}, Math. Comp. \textbf{38} (1982), no.~158, 517--529.

\bibitem{Sidi88}
\bysame, \emph{A user friendly extrapolation method for oscillatory
infinite
  integrals}, Math. Comp. \textbf{51} (1988), 249--266.

\bibitem{Wynn56}
P.~Wynn, \emph{On a device for computing the $e_m(s_n)$
tranformation},
  Mathematical Tables and Other Aids to Computation \textbf{10} (1956), 91--96.

\end{thebibliography}

\providecommand{\bysame}{\leavevmode\hbox
to3em{\hrulefill}\thinspace}
\providecommand{\MR}{\relax\ifhmode\unskip\space\fi MR }
\providecommand{\MRhref}[2]{%
  \href{http://www.ams.org/mathscinet-getitem?mr=#1}{#2}
} \providecommand{\href}[2]{#2}

\newpage

\begin{table}[htbp]
\begin{center}
\caption{Error of using formula (4), $\varepsilon_T$, in comparison
with the truncation error $I_T$ if formula (4) is not applied,
corresponding to Figure 1 in Example 1.}
\begin{tabular*}{0.5\textwidth}
{@{\extracolsep{\fill}}@{\hspace{0.02\textwidth}}ccc@{\hspace{0.02\textwidth}}}
\toprule  $\alpha$  & $\varepsilon_T$ & $I_T$ \\
 \midrule
 0.001 & $9.4\times 10^{-7}$ & 0.9391  \\
\midrule
 0.01 & $5.3\times 10^{-5}$ & 0.5334  \\
 \midrule
 0.1 & $1.8\times 10^{-5}$ & 0.0018  \\
 \midrule
 1 & $2.6\times 10^{-28}$ &  $2.6\times 10^{-28}$ \\
 \midrule
 10 & 0.0  & 0.0  \\
 \bottomrule
\end{tabular*}
\label{tab1}
\end{center}
\end{table}

\begin{table}[htbp]
\begin{center}
\caption{Error of using formula (4), $\varepsilon_T$, in comparison
with the truncation error $I_T$ if formula (4) is not applied,
corresponding to Figure 3 in Example 2.}
\begin{tabular*}{0.5\textwidth}
{@{\extracolsep{\fill}}@{\hspace{0.02\textwidth}}ccc@{\hspace{0.02\textwidth}}}
\toprule  $x$  & $\varepsilon_T$ & $I_T$ \\
 \midrule
 $4\pi$ & $1.0\times 10^{-3}$ & 0.2241  \\
\midrule
 $10\pi$& $1.1\times 10^{-4}$ & 0.1422  \\
 \midrule
 $20\pi$ & $1.9\times 10^{-5}$ & 0.1006  \\
 \midrule
 $50\pi$ & $1.9\times 10^{-6}$ &  0.0637 \\
 \midrule
 $100\pi$& $3.4\times 10^{-7}$  & 0.0318  \\
 \bottomrule
\end{tabular*}
\label{tab1b}
\end{center}
\end{table}

\begin{table}[htbp]
\begin{center}
\caption{Error of using formula (4), $\varepsilon_T$, in comparison
with the truncation error $I_T$ if formula (4) is not applied,
corresponding to Figure 5 in Example 3.}
\begin{tabular*}{0.5\textwidth}
{@{\extracolsep{\fill}}@{\hspace{0.02\textwidth}}ccc@{\hspace{0.02\textwidth}}}
\toprule  $x$  & $\varepsilon_T$ & $I_T$ \\
 \midrule
 $20\pi$ & $4.2\times 10^{-4}$ & 0.0105  \\
\midrule
 $40\pi$& $2.1\times 10^{-4}$ & 0.0053  \\
 \midrule
 $60\pi$ & $1.4\times 10^{-4}$ & 0.0035  \\
 \midrule
 $80\pi$ & $1.1\times 10^{-4}$ &  0.0026 \\
 \midrule
 $100\pi$& $8.4\times 10^{-5}$  & 0.0021  \\
 \bottomrule
\end{tabular*}
\label{tab1c}
\end{center}
\end{table}

\begin{figure}[htbp]
\centerline{\includegraphics{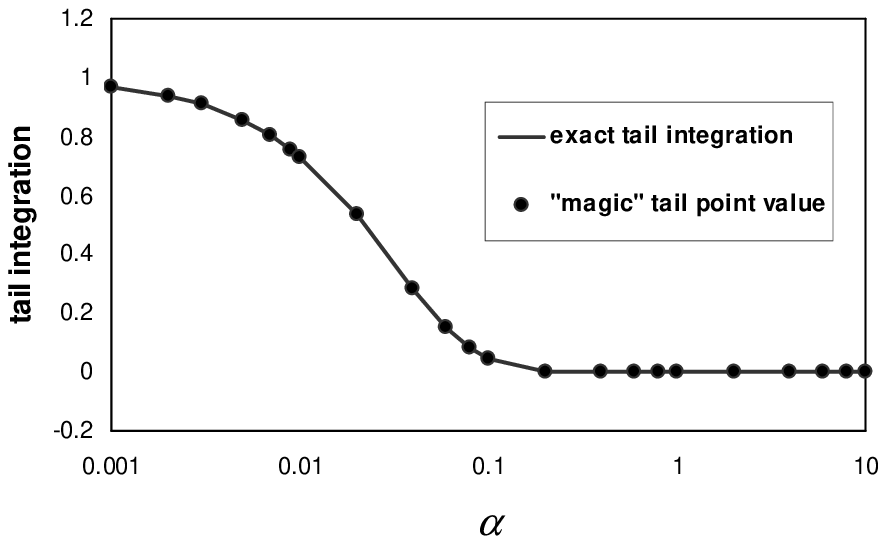}}  \caption{Comparison
between the exact tail integration $\int_{\pi N}^\infty {f(x)\sin
(x)dx} $ and simple one-point approximation $f(N\pi )$ from formula
(4), when $f(x) = e^{ - \alpha \,x}$ and $N = 10$.}\label{fig1}
\end{figure}

\bigskip

\begin{figure}[htbp]
\centerline{\includegraphics{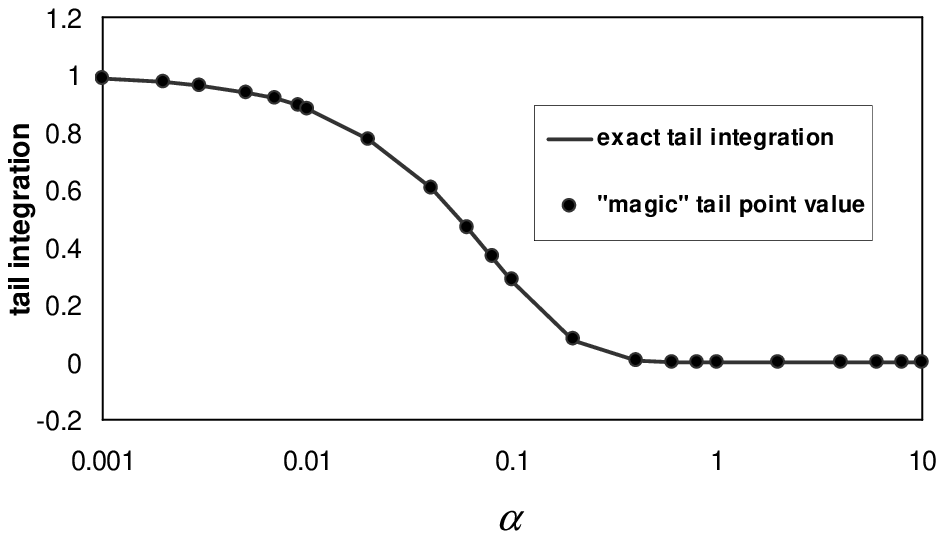}}  \caption{Comparison
between the exact tail integration $\int_{\pi N}^\infty {f(x)\sin
(x)dx} $ and simple one-point approximation $f(N\pi )$ from formula
(4), when $f(x) = e^{ - \alpha \,x}$ and $N = 4$.} \label{fig2}
\end{figure}

\bigskip

\begin{figure}[htbp]
\centerline{\includegraphics{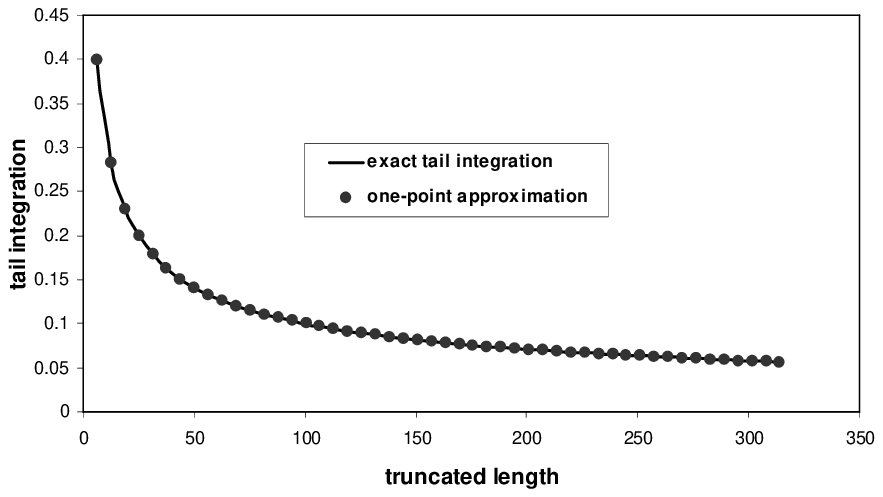}} \caption{Comparison
between exact tail integration $\int_{\pi N}^\infty {f(x)\sin (x)dx}
$ and the simple one point approximation (4), $f(N\pi )$, as
functions of truncated length $l_T = N\pi ,\;4 \le N \le 100$, when
$f(x) = 1 / \sqrt x $.} \label{fig3}
\end{figure}

\bigskip

\begin{figure}[htbp]
\centerline{\includegraphics{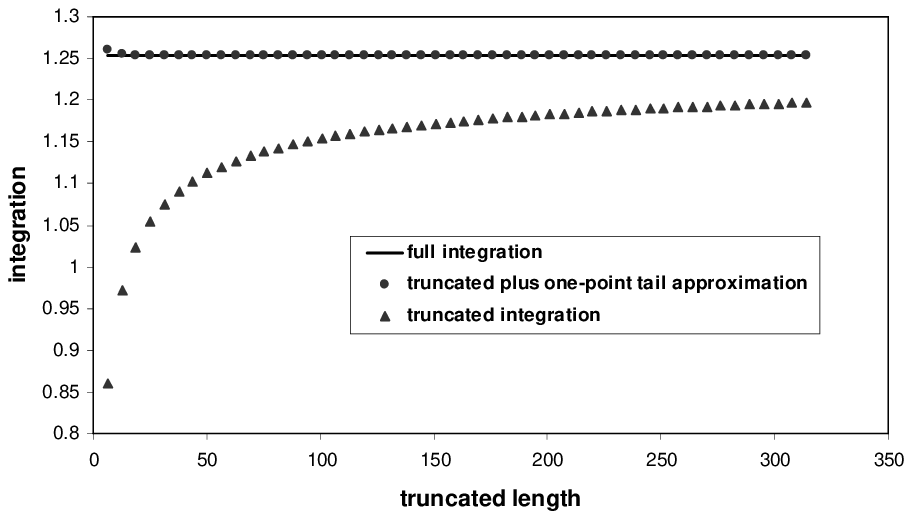}} \caption{Comparison
between truncated integration $\tilde {I}(N\pi ) = \int_0^{\pi N}
{f(x)\sin (x)dx} $ and the truncated integration plus the one-point
approximation of tail integration, $\tilde {I}(N\pi ) + f(N\pi )$,
as functions of the truncated length $l_T = N\pi ,\;4 \le N \le
100$, where $f(x) = 1/ \sqrt{x} $. The solid line represents the
exact value of the full integration without truncation error, $I_E =
\tilde {I}(\infty ) = \sqrt {\pi / 2} $.} \label{fig4}
\end{figure}

\bigskip

\begin{figure}[htbp]
\centerline{\includegraphics{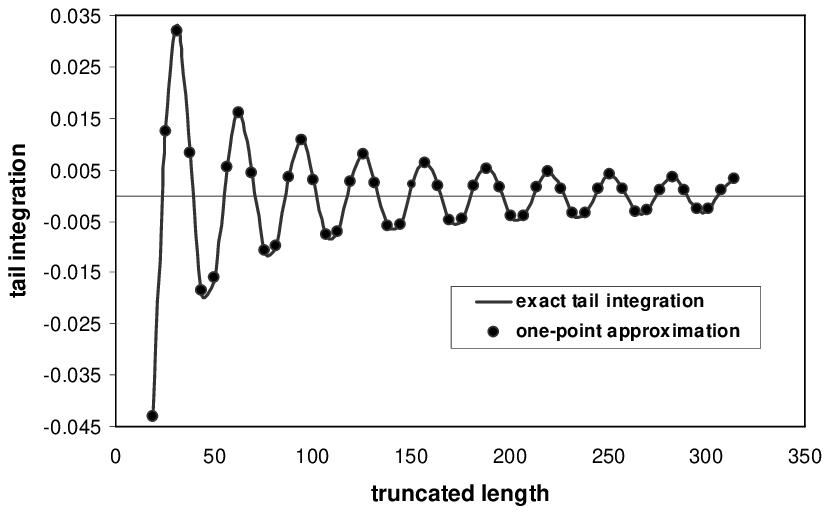}} \caption{Comparison
between exact tail integration $\int_{\pi N}^\infty {f(x)\sin (x)dx}
$ and the simple one-point approximation (4), $f(N\pi )$, as
functions of truncated length $l_T = N\pi ,\;6 \le N \le 100$, $f(x)
= \cos(\alpha x) / x $, $\alpha=0.2$.} \label{fig5}
\end{figure}

\bigskip

\begin{figure}[htbp]
\centerline{\includegraphics{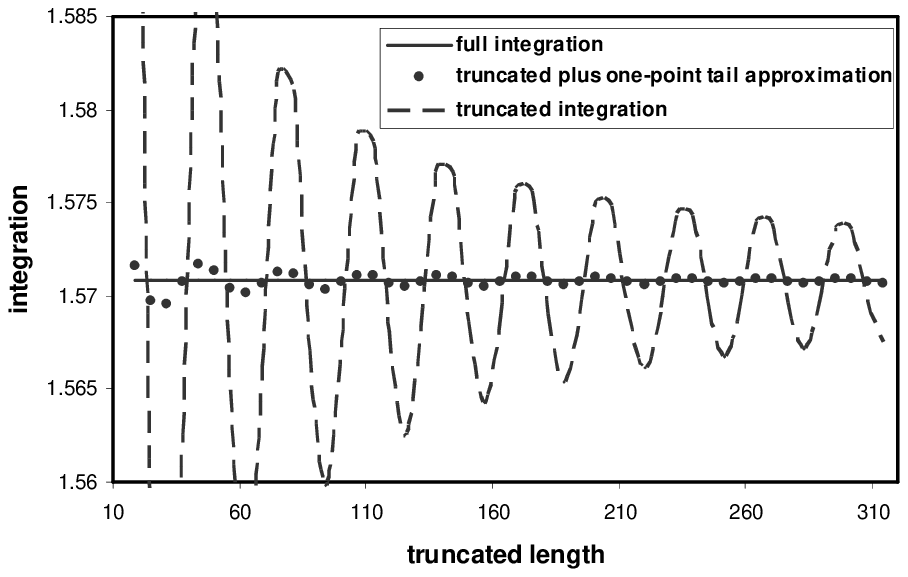}} \caption{Comparison
between truncated integration $\tilde {I}(N\pi ) = \int_0^{\pi N}
{f(x)\sin (x)dx} $ and the truncated integration plus the one-point
approximation of tail integration, $\tilde {I}(N\pi ) + f(N\pi )$,
as functions of the truncated length $l_T = N\pi ,\;6 \le N \le
100$, where $f(x) = \cos(\alpha x) / x $, $\alpha=0.2$. The solid
line represents the exact value of the full integration without
truncation error, $I_E = \tilde {I}(\infty ) = \pi / 2 $.}
\label{fig6}
\end{figure}

\end{document}